\documentclass[12pt]{article}

\usepackage[utf8]{inputenc}
\usepackage[utopia]{mathdesign}
\usepackage{latexsym}
\usepackage{amstext}
\usepackage{amsmath}
\usepackage{amsthm}
\usepackage{makeidx}
\usepackage{array}
\usepackage{graphicx}
\usepackage[all]{xy}

\newtheorem{thm}{Theorem}[section]

\newtheorem{prop}[thm]{Proposition}
\newtheorem{cor}[thm]{Corollary}
\theoremstyle{definition}
\newtheorem{defn}{Definition}[section]

\title{Artin algebras of finite type and finite categories of $\Delta$-good modules}

\author{Danilo D. da Silva\\
\small Universidade Federal de Sergipe - UFS\\[-0.8ex]
\small Departamento de Matematica - DMA\\[-0.8ex]
\small Sao Cristovao, SE, Brazil\\
\small \texttt{ddsilva@ufs.br}\\
}

\date{{Feb 13, 2015}}

\begin{document}
\maketitle

\begin{abstract}
We give an alternative proof to the fact that if the square of the infinite radical of the module category of an Artin algebra is equal to zero then the algebra is of finite type by making use of the theory of postprojective and preinjective partitions. Further, we use this new approach in order to get a characterization of finite subcategories of $\Delta$-good modules of a quasi-hereditary algebra  in terms of depth of morphisms similar to a recently obtained characterization of Artin algebras of finite type.  
\end{abstract}

%\subjclass{16G70, 16G60, 16G10}

\section{Introduction}

Let $A$ be an Artin algebra. In \cite{coelhoetall} it is shown that if the square of the infinite radical of the category ${\rm mod}A$ is equal to zero, $({\rm rad}_A^{\infty})^2=0$, then $A$ is of finite type. This result was used then in \cite{chaioliu} to obtain a characterization of Artin algebras of finite type in terms of depths of morphisms which, in particular, improved the result proved in \cite{coelhoetall}.

On the other hand, it is shown in \cite{auslander} that the category ${\rm ind}A$ has postprojective and preinjective partitions (in \cite{auslander} the postprojective partition was called then preprojective partition). Further, in \cite{skowronskismalo} it is proved that being an algebra of infinite type there always exists an indecomposable module which is neither postprojective nor preinjective. In this article we use the theory of postprojective and preinjective partitions and the result proved in \cite{skowronskismalo} to obtain an alternative proof to the fact that if $({\rm rad}_{A}^{\infty})^2=0$ then $A$ is of finite type (Corollary \ref{corolario}). 

In the study of a quasi-hereditary algebra $A$ is of big importance to find conditions of finiteness for the category ${\cal F}(\Delta)$ of $\Delta$-good modules which consists of $A$-modules which have a filtration by standard modules. Recently, it has been given attention to which sufficient finiteness conditions for ${\rm mod}A$ could also be proved true for ${\cal F}(\Delta)$. For instance, in \cite{deng} it was proved for ${\cal F}(\Delta)$ a similar finiteness condition to the one proved in \cite{happel}, that is, if ${\cal F}(\Delta)$ does not have short cycles then ${\cal F}(\Delta)$ is of finite type. In \cite{zeng}, it is proved that a finite-dimensional quasi-hereditary algebra over an infinite perfect field satisfies Brauer-Thrall II and, as consequence, that if each module in ${\cal F}(\Delta)$ is either postprojective or preinjective then ${\cal F}(\Delta)$ is finite.

In this paper we intend to use the result proved in \cite{zeng} mentioned above to find similar results to the ones obtained in \cite{coelhoetall} and \cite{chaioliu} for ${\cal F}(\Delta)$. More explicitly, if ${\rm rad}_{\Delta}$ is the radical of the category ${\cal F}(\Delta)$  we prove that if $({\rm rad}_{\Delta}^{\infty})^2=0$ then ${\cal F}(\Delta)$ is finite (Theorem \ref{thmdan}) and moreover we also prove a characterization of finiteness for ${\cal F}(\Delta)$ similar to the one found in \cite{chaioliu} (Theorem \ref{thmcharac}).

We finish the paper interested in finding a bound to the number of levels of the postprojective partition (similarly, the preinjective partition) whenever ${\cal F}(\Delta)$ is finite (Propositions \ref{boundpost} and \ref{boundpre}).

\section{Preliminaries}

Let $A$ be an Artin algebra and let ${\rm mod}A$ be the category of finitely generated right $A$-modules. By a subcategory ${\mathfrak C}$ of ${\rm mod}A$ we always mean a full subcategory of ${\rm mod}A$ closed under direct summands. We denote the subcategory of ${\rm mod}A$ of indecomposable $A$-modules by ${\rm ind}A$ and the subcategory of indecomposable $A$-modules in ${\mathfrak C}$ by ${\rm ind}{\mathfrak C}$. The category ${\rm add}{\mathfrak C}$ is the subcategory of ${\rm mod}A$ consisting of all $A$-modules isomorphic to summands of finite sums of modules in ${\mathfrak C}$.

We say a subategory ${\mathfrak C}$ of ${\rm mod}A$ is {\it covariantly finite} in ${\rm mod}A$ if given an $A$-module $D$ in ${\rm mod}A$ there exists a morphism $f: D \rightarrow C$ with $C \in {\rm add}{\mathfrak C}$ such that for every $X \in {\mathfrak C}$ we have that ${\rm Hom}(f,X):{\rm Hom}(C,X) \rightarrow {\rm Hom}(D,X)$ is surjective. Dually, we have the concept of a {\it contravariantly finite} subcategory in ${\rm mod}A$. Finally, we say the subcategory ${\mathfrak C}$ is {\it functorially finite} in ${\rm mod}A$ whenever it is both covariantly finite and contravariantly finite in ${\rm mod}A$. 

Let ${\mathfrak C}$ be a subcategory of ${\rm mod}A$. Given two $A$-modules $M$ and $N$ in ${\mathfrak C}$, we define ${\rm Hom}_{{\mathfrak C}}(M,N)={\rm Hom}_A(M,N)$ and we consider ${\rm rad}_{{\mathfrak C}}(M,N)$ to be the set of morphisms $f \in {\rm Hom}_{{\mathfrak C}}(M,N)$ such that for any $X \in {\rm add}{\mathfrak C}$ indecomposable and for all morphisms $g: X \rightarrow M$ and $h: N \rightarrow X$ we have that $hfg$ is not an isomorphism. Given $m > 1$, we define ${\rm rad}_{{\mathfrak C}}^{m}(M,N) \subseteq {\rm rad}_{{\mathfrak C}}(M,N)$ as the subgroup consisting of all finite sums of morphisms of the form

$$M=M_0 \stackrel{h_1}{\rightarrow} M_1 \stackrel{h_2}{\rightarrow}M_2 \rightarrow \cdots \rightarrow M_{m-1} \stackrel{h_{m}}{\rightarrow} M_{m}=N$$

where $h_j \in {\rm rad}_{{\mathfrak C}}(M_{j-1},M_j)$. Thus we get the chain:
$${\rm Hom}_{{\mathfrak C}}(M,N) \supseteq {\rm rad}_{{\mathfrak C}}(M,N) \supseteq {\rm rad}_{{\mathfrak C}}^2(M,N) \supseteq \cdots \supseteq {\rm rad}_{{\mathfrak C}}^n(M,N) \supseteq \cdots$$

Finally, we set ${\rm rad}_{{\mathfrak C}}^{\infty}(M,N)= \bigcap_{n \geq 1}{\rm rad}_{{\mathfrak C}}^n(M,N)$. When ${\mathfrak C}={\rm mod}A$, ${\rm rad}_{\mathfrak C}$ is the usual radical of the category ${\rm mod}A$ denoted by ${\rm rad}_A$.

\begin{defn}
Let ${\mathfrak C}$ be a subcategory of ${\rm mod}A$ and $f:X \rightarrow Y$ a morphism in ${\rm rad}_{{\mathfrak C}}(X,Y)$. Let $n$ be a positive integer. We say the {\it depth} of $f$ relative to ${\mathfrak C}$ is equal to $n$ and we denote ${\rm dp}_{{\mathfrak C}}(f) = n$ if $f \in {\rm rad}_{{\mathfrak C}}^{n}(X,Y) \backslash {\rm rad}_{{\mathfrak C}}^{n+1}(X,Y)$. If there is no such $n$ we say the depth of $f$ relative to ${\mathfrak C}$ is infinite and we denote ${\rm dp}_{{\mathfrak C}}(f) = \infty$.
\end{defn}
If ${\mathfrak C}={\rm mod}A$ then we have the usual notion of depth of a morphism $f$ introduced in \cite{chaioliu} denoted there by ${\rm dp}(f)$.

Suppose ${\mathfrak C}$ is a subcategory of ${\rm mod}A$ and $\cal{X}$ a subcategory of ${\mathfrak C}$. We denote by ${\mathfrak C}_{\cal{X}}$ the subcategory of ${\mathfrak C}$ consisting of the objects in ${\mathfrak C}$ with no summands in $\cal{X}$. We say that a module $N$ in ${\mathfrak C}$ is a {\it splitting projective} in ${\mathfrak C}$ if each epimorphism $M \rightarrow N$ with $M \in {\rm add}{\mathfrak C}$ splits. We denote by ${\bf P_0({\mathfrak C})}$ the subcategory of ${\rm ind}{\mathfrak C}$ consisting of the indecomposable splitting projectives in ${\mathfrak C}$. We define ${\bf P_1}({\mathfrak C})={\bf P_0}({\mathfrak C}_{{\bf P_0({\mathfrak C})}})$ and, by induction, ${\bf P_k}({\mathfrak C})={\bf P_0({\mathfrak C}_{{\bf P_0}({\mathfrak C}) \cup \cdots \cup {\bf P_{k-1}}({\mathfrak C})})}$. Finally, we denote the subcategory $\bigcup_{i<{\infty}}{\bf P_i}({\mathfrak C})$ of ${\rm ind}{\mathfrak C}$ by ${\bf P}({\mathfrak C})$ and we set ${\bf P_{\infty}}({\mathfrak C})={\rm ind}({\mathfrak C}_{{\bf P}({\mathfrak C})})$, the subcategory whose objects are the indecomposable modules in $\rm{add}{\mathfrak C}_{{\bf P}({\mathfrak C})}$. 

A {\it cover} for a subcategory ${\mathfrak C}$ of ${\rm ind}A$ is a subcategory ${\mathfrak D}$ of ${\mathfrak C}$ such that for each $X$ in ${\mathfrak C}$ there is a surjective morphism $f:Y \rightarrow X$ with $Y$ in ${\rm add}{\mathfrak D}$. It is called a minimal cover if no proper subcategory of ${\mathfrak D}$ is a cover for ${\mathfrak C}$.

We say the collection $\{ {\bf P_i}({\mathfrak C})\}_{i=0, \cdots , {\infty}}$ as above is a {\it postprojective partition} of ${\mathfrak C}$ if ${\bf P_i}({\mathfrak C})$ is a finite cover for ${\mathfrak C}_{{\bf P_0}({\mathfrak C}) \cup \cdots \cup {\bf P_{i-1}}({\mathfrak C})}$ for each $i< {\infty}$. If ${\mathfrak C}$ has a postprojective partition we say that $M \in {\mathfrak C}$ is a {\it postprojective module} (former preprojective module in \cite{auslander}) if every indecomposable summand of $M$ is in ${\bf P}({\mathfrak C})$. 

If an indecomposable $X \in {\mathfrak C}$ is in ${\bf
P_n}({\mathfrak C})$, $0 \leq n < {\infty}$, then we say that $X$ is postprojective of {\it level} $n$.
Finallly, we denote by ${\bf P^m}({\mathfrak C})$ the subcategory ${\bf P_0}({\mathfrak C}) \cup \cdots \cup {\bf P_{m}}({\mathfrak C})$.

Dually, we have the concepts of {\it splitting injectives}, {\it cocovers}, {\it preinjective partition} and {\it preinjective modules}. For a more detailed account on the theory we refer the reader to \cite{auslander}.

I turns out (see \cite{auslander}) that a covariantly finite subcategory in ${\rm mod}A$ has always a uniquely determined postprojective partition and that a contravariantly finite subcategory in ${\rm mod}A$ has always a uniquely determined preinjective partition.

We say that a subcategory of ${\rm mod}A$ is {\it resolving} if it contains all the projective indecomposable $A$-modules and if it is closed under extensions and also under kernels of epimorphisms. Dually, we can define a {\it coresolving} subcategory of ${\rm mod}A$.

The following proposition shall be used very often in this paper.

\begin{prop}
\label{propdan}
Let $A$ be an Artin algebra and let ${\mathfrak C}$ be a resolving subcategory covariantly finite in ${\rm mod}A$. Given $0<i \leq {\infty}$, we have ${\rm Hom}_{{\mathfrak C}}(M,N)={\rm rad}_{{\mathfrak C}}^i(M,N)$ $\forall M \in \bf{P_0({\mathfrak C})}$, $\forall N \in \bf{P_i({\mathfrak C})}$.

\begin{proof}
We first prove, by induction on $i$, that given $1\leq i<{\infty}$ we have\break ${\rm Hom}_{{\mathfrak C}}(M,N)={\rm rad}_{{\mathfrak C}}^i(M,N)$ $\forall M \in \bf{P_0}({\mathfrak C})$, $\forall N \in \bf{P_i}({\mathfrak C})$. As $\bf{P_0}({\mathfrak C}) \cap \bf{P_1}({\mathfrak C})=\emptyset$ it is clear, for all $M \in \bf{P_0}({\mathfrak C})$ and $N \in \bf{P_1}({\mathfrak C})$, that ${\rm Hom}_{{\mathfrak C}}(M,N)={\rm rad}_A(M,N) \subseteq {\rm rad}_{{\mathfrak C}}(M,N)$. Hence ${\rm Hom}_{{\mathfrak C}}(M,N)={\rm rad}_{{\mathfrak C}}^1(M,N)$ $\forall M \in \bf{P_0({\mathfrak C})}$, $\forall N \in \bf{P_1({\mathfrak C})}$. Suppose we have the equality for $i-1$ and take $f : M_0 \rightarrow M_i$ with $M_0 \in \bf{P_0}({\mathfrak C})$ and $M_i \in \bf{P_i}({\mathfrak C})$. We know $\bf{P_{i-1}}({\mathfrak C})$ is a cover for ${\mathfrak C}_{{\bf P^{i-2}}({\mathfrak C})}$ so we get an epimorphism $\displaystyle\bigoplus_{l=1}^{n}M_{{i-1},l} \stackrel {[h_1, \cdots, h_n]}{\longrightarrow} M_i$ with $M_{{i-1},l} \in \bf{P_{i-1}}({\mathfrak C})$, $\forall l \in \{1, \cdots, n \}$. Since ${\mathfrak C}$ is a resolving subcategory, $\bf{P_0({\mathfrak C})}$ consists of all the indecomposable projective $A$-modules so we get the lifting:

$$\xymatrix{
                                                                                                                &      &        M_0  \ar[lld]_-{[g_1, \cdots, g_n]^t} \ar[d]^{f}                    &              \\
\displaystyle\oplus_{l=1}^{n}M_{{i-1},l}  \ar[rr]_-{[h_1, \cdots, h_n]}     &        &      M_i     \ar[r]                   & 0             \\
}$$

As $\bf{P_{i-1}}({\mathfrak C}) \cap \bf{P_i}({\mathfrak C})=\emptyset$, we know that $h_l : M_{{i-1},l} \rightarrow M_i$ is in ${\rm rad}_{{\mathfrak C}}(M_{i-1,l},M_i)$ for every $l \in \{1, \cdots, n \}$ and we also know, by hypothesis, that $g_l \in {\rm rad}_{{\mathfrak C}}^{i-1}(M_0,M_{{i-1},l}) $ for every $l \in \{1, \cdots, n \}$. Hence $ f=\displaystyle\sum_{l=1}^n h_l g_l \in {\rm rad}^{i}_{\mathfrak C}(M_0,M_{i}) $ and ${\rm Hom}_{{\mathfrak C}}(M,N)={\rm rad}_{{\mathfrak C}}^i(M,N)$ $\forall M \in \bf{P_0}({\mathfrak C})$, $\forall N \in \bf{P_i}({\mathfrak C})$.
 
Consider now $f: M_0 \rightarrow M_{\infty}$ with $M_0\in \bf{P_0}({\mathfrak C})$ and $M_{\infty}\in \bf{P_{\infty}}({\mathfrak C})$. Given $0<i<{\infty}$, since ${\bf P_i}({\mathfrak C})$ is a cover for ${\mathfrak C}_{{\bf P^{i-1}}({\mathfrak C})}$ there exists an epimorphism $\displaystyle\bigoplus_{k=1}^{n}M_{{i},l} \stackrel {[h_1, \cdots, h_n]}{\longrightarrow} M_{\infty}$ with $M_{{i},l} \in \bf{P_{i}}({\mathfrak C})$, $\forall l \in \{1, \cdots, n \}$. Again we get the lifting:
%\vspace{1 cm}\\

$$\xymatrix{
                                                                                                                &      &        M_0  \ar[lld]_-{[g_1, \cdots, g_n]^t} \ar[d]^{f}                    &              \\
\displaystyle\oplus_{l=1}^{n}M_{{i},l}  \ar[rr]_-{[h_1, \cdots, h_n]}     &        &      M_{\infty}     \ar[r]                   & 0             \\
}$$

We know already that $g_l \in {\rm rad}_{{\mathfrak C}}^i(M_0,M_{i,l})$ $\forall l \in \{1, \cdots, n \}$ and hence $f \in {\rm rad}_{{\mathfrak C}}^{i+1}(M_0,M_{\infty})$. As $i$ was taken arbitrarily, we get $f \in {\rm rad}_{{\mathfrak C}}^{\infty}(M_0,M_{\infty})$.

\end{proof}

\end{prop}

Duallly, if we consider a coresolving contravariantly finite subcategory of ${\rm mod}A$ and its preinjective partition ${\bf I_0}({\mathfrak C}), {\bf I_1}({\mathfrak C}), \cdots, {\bf I_{\infty}}({\mathfrak C})$ then given $0 < j \leq {\infty}$ we have ${\rm Hom}_{{\mathfrak C}}(M,N)={\rm rad}_{{\mathfrak C}}^j(M,N)$ $\forall M \in \bf{I_j}({\mathfrak C})$, $\forall N \in \bf{I_0}({\mathfrak C})$.\\

\section{Artin algebras of finite type}

In this section, let $A$ be an Artin algebra unless otherwise stated. As mentioned before, the goal of this section is to use the theory of postprojective and preinjective partitions to prove a known result obtained in \cite{coelhoetall}, namely if ${\rm rad}^2_A=0$ then $A$ is of finite type. This new approach will also be shown to be useful to prove a similar result for subcategories ${\cal F}(\Delta)$ of ${\rm mod}A$ where $A$ is quasi-hereditary.  

We denote by ${\bf P}_0, {\bf P}_1, \cdots, {\bf P}_{\infty}$ and ${\bf I}_0, {\bf I}_1, \cdots, {\bf I}_{\infty}$ the postprojective and the preinjective partitions of ${\rm ind}A$, respectively.

We begin by recalling the next theorem which is key to our approach.

\begin{thm}

{\rm \cite{skowronskismalo}}
\label{skowronskismalo}
If \, ${\bf P_{\infty} \cap {\bf I_{\infty}}=\emptyset}$ then $A$ is of finite type.
\end{thm}

Following \cite{chaioliu}, given a simple $A$-module $S$ we fix a projective cover $\pi_S:P_S \rightarrow S$ and an injective envelope $\iota_S : S \rightarrow I_S$. Moreover, we set $\theta_S=\iota_S \pi_S$. 

The result proved in \cite{coelhoetall} mentioned above was used in \cite{chaioliu} to prove the equivalence between itens a) and d) below. 

\begin{prop}
\label{chaioliu}
\cite{chaioliu}
The following are equivalent:
\begin{itemize}

\item[a)] $A$ is of finite type.
\item[b)] The depth of $\pi_S$ is finite for every simple $A$-module $S$.
\item[c)] The depth of $\iota_S$ is finite for every simple $A$-module $S$.
\item[d)] $\theta_S \not\in ({\rm rad_A}^{\infty})^2$ for every simple $A$-module $S$.
\end{itemize}

\end{prop}

Next we give another proof to the equivalence between itens a) and d) without using the result proved in \cite{coelhoetall}. In the next section we shall prove a similar result to Proposition \ref{chaioliu} for subcategories of $\Delta$-good modules.

\begin{prop}
\label{novaprop}
The following are equivalent:
\begin{description}
\item [a)] A is of finite type.
\item [b)] $\theta_S \not\in ({\rm rad_A}^{\infty})^2$ for every simple $A$-module $S$.
\item [c)] $({\rm rad}_A^{\infty})^2(P_S,I_S)=0$ for every simple $A$-module $S$.
\end{description}
\begin{proof}

We observe first that ${\rm mod}A$ is resolving and coresolving and that it is also functorially finite in itself. 

Clearly item a) implies item c). We also have that item c) implies item b) since $\theta_S \neq 0$ for every simple $A$-module $S$. So we need to prove that item b) implies item a).

Suppose by contradiction that $A$ is of infinite type. Then we have ${\bf P}_{\infty} \cap {\bf I}_{\infty} \neq \emptyset$, by Theorem \ref{skowronskismalo}. Take $M \in {\bf P}_{\infty} \cap {\bf I}_{\infty}$ and consider the inclusion $g : S \hookrightarrow M$, where $S$ is a simple $A$-module summand of the socle of $M$. Since $M \in {\bf I_{\infty}}$ we get $S \in {\bf I_{\infty}}$ (it could happen $M=S$). Observe the diagram:

$$\xymatrix{
            P_S \ar[r]^{\pi_S}                 &  S \ar[r]^{\iota_S} \ar[d]^{g}   &  I_S          \\
                                                        &   M     \ar[ru]^{h}                     &                  \\
}$$

where $h$ is obtained from the fact that $I_S$ is injective. 

We have then $\theta_S=\iota_S \pi_S= h(g\pi_S)$ and, hence, $\theta \in ({\rm rad_A}^{\infty})^2$ since $h \in \break {\rm Hom}_A(M,I_S)={\rm rad_A}^{\infty}(M,I_S)$, by the dual of Proposition \ref{propdan}, and $g\pi_S \in \break {\rm Hom}_A(P_S,M)={\rm rad_A}^{\infty}(P_S,M)$, by Proposition \ref{propdan} itself. Observe that if $P_S=S$ then $S \neq M$, $\pi_S=1_S$, $\theta_S=\iota_S$ and $g \in {\rm rad_A}^{\infty}(S,M)$ since $S$ is projective and $M \in {\bf P_{\infty}}$ so once more we get $\theta_S \in ({\rm rad_A}^{\infty})^2$.

\end{proof}
\end{prop}

\begin{cor}
\label{corolario}
$A$ is of finite type if and only if $({\rm rad}^{\infty}_A)^2=0$.
\begin{proof}
It follows as a consequence of the equivalence between itens $c)$ and $a)$ in Proposition \ref{novaprop}.
\end{proof}
\end{cor}

\section{Partitions of a category ${\cal F}(\Delta)$ of $\Delta$-good modules}

From this point on, we assume that $A$ is quasi-hereditary with a fixed ordering $S(1), S(2), \cdots, S(n)$ of the isomorphism classes of the simple $A$-modules. We also assume that $\Delta(1), \cdots, \Delta(n)$ are the corresponding standard $A$-modules and $T(1), T(2), \cdots, T(n)$ are the characteristic $A$-modules (see \cite{ringel}). Let ${\cal F}(\Delta)$ be the subcategory of $\Delta$-good modules of ${\rm mod}A$ which, by definition, consists of $A$-modules having a $\Delta$-good filtration, that is, a filtration by the standard $A$-modules. 

It is proved in \cite{ringel} that ${\cal F}(\Delta)$ is functorially finite in ${\rm mod}A$, closed under extensions and closed under direct summands. In particular, ${\cal F}(\Delta)$ has relative almost split sequences and it admits both postprojective and preinjective partitions, denoted here by ${\bf P_0}(\Delta), {\bf P_1}(\Delta), \cdots , {\bf P_{\infty}}(\Delta)$ and ${\bf I_0}(\Delta), {\bf I_1}(\Delta), \cdots , {\bf I_{\infty}}(\Delta)$, respectively. Further, it was also proved in \cite{ringel} that ${\cal F}(\Delta)$ is a resolving subcategory of ${\rm mod}A$. 

Moreover, in \cite{zeng} it was proved the following:

\begin{thm}
\label{teozeng}
{\rm \cite{zeng}}
Let $A$ be a finite dimensional quasi hereditary algebra over a infinite perfect field $k$. If ${\cal F}(\Delta)$ is infinite then ${\bf P}_{\infty}(\Delta) \cap {\bf I}_{\infty}(\Delta) \neq \emptyset$.
\end{thm}

As mentioned before, we shall prove for ${\cal F}(\Delta)$ a similar result to the one proved in \cite{coelhoetall} under the hypothesis that $A$ is a finite dimensional quasi-hereditary algebra over an infinite perfect field since we shall need Theorem \ref{teozeng} in the proof.

The next proposition follows directly from the definitions and it appeared first in \cite{zeng}.

\begin{prop}
${\bf P_0}(\Delta)$ consists of all indecomposable projective modules and ${\bf I_0}(\Delta)$ of all indecomposable {\rm Ext}-injective modules, that is, the characteristic modules $T(1), \cdots, T(n)$.

\end{prop}

Denoting ${\rm Hom}_{\Delta}(M,N)={\rm Hom}_{{\cal F}(\Delta)}(M,N)$, for every $0<i \leq {\infty}$, since ${\cal F}(\Delta)$ is a resolving subcategory we have, by Proposition \ref{propdan}, ${\rm Hom}_{\Delta}(M,N)=\break{\rm rad}_{\Delta}^i(M,N)$ $\forall M \in {\bf P_0}(\Delta)$, $\forall N \in {\bf P_i}(\Delta)$. We do not get the dual version this time because $T(1), \cdots, T(n)$ are not necessarily injective (they are Ext-injective modules in a subcategory that is not necessarily closed under cokernels), so we lose the lifting property.

Now we shall prove one of the main results of the paper, namely if $({\rm rad}_{\Delta}^{\infty})^2=0$ then ${\cal F}(\Delta)$ is finite. But before this we need the following proposition:

\begin{prop}
\label{propreview}
If $M$ is indecomposable in ${\bf I}_{\infty}(\Delta)$, then there exists a module $T \in\, {\bf I}_{0}(\Delta)$ with ${\rm rad}^{\infty}_{\Delta}(M,T) \neq 0$.
\begin{proof}
For each $m \in \mathbb{N}$ there exists a chain of monomorphisms 

 $$ \rm{(} \ast \rm{)} \;\;\;\;\;  M \rightarrow Z_m \rightarrow Z_{m-1} \rightarrow \cdots \rightarrow Z_1 \rightarrow Z_0$$

with $Z_i \in {\rm add}{\bf I}_i(\Delta)$. Since ${\bf I}_i(\Delta) \cap {\bf I}_{j}(\Delta)=\emptyset$ for $i \neq j$, this chain of monomorphisms is in ${\rm rad}_{\Delta}^{m+1}(M,Z_0)$. Consequently there exists a module $T(i_m) \in {\bf I}_{0}(\Delta)$ with ${\rm rad}_{\Delta}^{m+1}(M,T(i_m)) \neq 0$. Since ${\bf I}_0(\Delta)$ is a finite set and Hom-spaces are finite dimensional, there exists a module $T \in {\bf I}_0(\Delta)$ with ${\rm rad}^{\infty}_{\Delta}(M,T) \neq 0$.

\end{proof}
\end{prop}

\begin{thm}
\label{thmdan}
Let $A$ be a quasi-hereditary algebra over an infinite perfect field $k$. If $({\rm rad}_{\Delta}^{\infty})^2=0$ then ${\cal F}(\Delta)$ is finite.

\begin{proof}
Suppose ${\cal F}(\Delta)$ is infinite. We know there exist a postprojective partition ${\bf P}_0(\Delta) \cdots {\bf P}_{\infty}(\Delta)$ and a preinjective partition ${\bf I}_0(\Delta) \cdots {\bf I}_{\infty}(\Delta)$ of ${\cal F}(\Delta)$ such that ${\bf P}_{\infty}(\Delta) \cap {\bf I}_{\infty}(\Delta) \neq \emptyset$, by Theorem \ref{teozeng}. Thus we can take $M \in {\bf P}_{\infty}(\Delta) \cap {\bf I}_{\infty}(\Delta)$. Since all the projectives $A$-modules are in ${\cal F}(\Delta)$ we consider an epimorphism $f:P \rightarrow M$ where $P$ is projective so we have $f \in {\rm rad}_{\Delta}^{\infty}(P,M)$, by Proposition \ref{propdan}. By Proposition \ref{propreview}, there exist a module $T \in {\bf I}_0(\Delta)$ and a morphism $g \neq 0$ in ${\rm rad}^{\infty}_{\Delta}(M,T)$. Therefore, $gf \in ({\rm rad}_{\Delta}^{\infty})^2$ and $gf \neq 0$ since $f$ is an epimorphism.

\end{proof}

\end{thm}

Below, the standard modules $\Delta(i)$ shall play similar role to the simple modules in order to obtain a characterization equivalent to the one in Proposition \ref{chaioliu}, but this time for the subcategory ${\cal F}(\Delta)$. The morphisms $\pi(i): P(i) \rightarrow \Delta(i)$ for all $i \in \{1, \cdots , n \}$ given by the definition of ${\cal F}(\Delta)$ shall play the role of the projective coverings of the simple modules and the morphisms $\beta(i): \Delta(i) \rightarrow T(i)$ for every $i \in \{1, \cdots , n \}$ given by Proposition 2 in \cite{ringel} shall play the role of the injective envelopes of the simple modules.

From this point on we fix morphisms $\beta(i): \Delta(i) \rightarrow T(i)$ for every $i \in \{1, \cdots , n \}$ as given by Proposition 2 in \cite{ringel} and we denote ${\rm dp}_{\Delta}(f)={\rm dp}_{{\cal F}(\Delta)}(f)$.

\begin{thm}
\label{thmcharac}
Let $A$ be a finite dimensional quasi-hereditary algebra over a field $k$. The following are equivalent:

\begin{itemize}
\item[a)] ${\cal F}(\Delta)$ is finite.
\item[b)] The depth of $\pi(i)$ related to ${\cal F}(\Delta)$ is finite for all $i \in \{1, \cdots , n \}$.
\item[c)] The depth of $\beta(i)$ related to ${\cal F}(\Delta)$ is finite for all $i \in \{1, \cdots , n \}$.
\end{itemize}
If additionally, $k$ is an infinite perfect field, then the itens above are equivalent to:
\begin{itemize}
\item[d)] $({\rm rad}_{\Delta}^{\infty})^2=0$
\end{itemize}
\begin{proof}

Since b)-d) trivially follow from a) it is enough to show the reverse implications.

$b) \Rightarrow a)$ Suppose by contradiction that ${\cal F}(\Delta)$ is infinite. 

Then there exists  $i \in \{1, \cdots, n \}$ such that $\Delta(i) \in {\bf P_{\infty}}(\Delta)$: in fact, otherwise we could take $m>0$ such that $\Delta(i) \in {\bf P^m}(\Delta)$ for all $i \in \{1, \cdots, n \}$. Thus taking arbitrarily $M \in {\bf P_{m+1}}(\Delta)$ we get a filtration $0=M_0 \subset \cdots \subset M_s \subset M_{s+1}=M$ which give us an epimorphism $M \rightarrow M_{s+1}/M_s = \Delta(j)$ for some $j \in \{1, \cdots, n \}$. This epimorphism splits since $\Delta(j) \in {\bf P^m}(\Delta)$ and this is a contradiction because $M$ is indecomposable.  

Thus if we take $i \in \{1, \cdots, n \}$ with $\Delta(i) \in {\bf P_{\infty}}(\Delta)$ then by Proposition \ref{propdan} the epimorphism $\pi(i): P(i) \rightarrow \Delta(i)$ has infinite depth related to ${\cal F}(\Delta)$ which is a contradiction.

$c) \Rightarrow a)$ Again suppose, by contradiction, that ${\cal F}(\Delta)$ is infinite. 

In a similar way we can show that there exists $i \in \{1, \cdots, n \}$ such that $\Delta(i) \in {\bf I_{\infty}}(\Delta)$. By Proposition 2 in \cite{ringel}, we know there exists an exact sequence 

$$0 \rightarrow \Delta(i) \stackrel{\beta(i)} {\rightarrow} T(i) \rightarrow X(i) \rightarrow 0$$

with  $X(i) \in {\cal F}(\Delta)$. Moreover, since $\Delta(i) \in {\bf I_{\infty}}(\Delta)$ there exist - from ($\ast$) in the proof of Proposition \ref{propreview}, a monomorphism $g: \Delta(i) \rightarrow T$ such that $T \in {\bf I}_0(\Delta)$ and $g \in {\rm rad}_{\Delta}^{\infty}(\Delta(i), T)$. Then consider the diagram with exact sequences:
$$\xymatrix{
      0 \ar[r]      & \Delta(i) \ar[r]^{\beta(i)} \ar[d]^{g}         &  T(i) \ar[r] \ar[d]^{g'}           &   X(i) \ar[r] \ar[d]^{1_{X(i)}}                         &   0      \\
      0 \ar[r]      & T \ar[r]^{h}                              &  W \ar[r]                              &   X(i) \ar[r]                                &   0      \\
}$$
where $W$ is the push-out of $T \stackrel{g}{\leftarrow} \Delta(i) \stackrel{\beta(i)}{\rightarrow} T(i)$. Since ${\cal F}(\Delta)$ is closed under extensions we get $W \in {\cal F}(\Delta)$. Moreover, we see from the diagram that if $g$ is a monomorphism then $g'$ is also a monomorphism. Thus, we have that $g'$ splits since $T(i) \in {\bf I_0}(\Delta)$ and $W \in {\cal F}(\Delta)$. If we take $h': W \rightarrow T(i)$ such that $h'g'=1_{T(i)}$ then $\beta(i)=h'hg$ which implies $\beta(i) \in {\rm rad}_{\Delta}^{\infty}(\Delta(i), T(i))$ since $g \in {\rm rad}_{\Delta}^{\infty}(\Delta(i), T)$. This contradicts the hypothesis of item c).

$d) \Rightarrow a)$ Theorem \ref{thmdan}.

\end{proof}
\end{thm}

We are now interested in finding a bound to the number of levels of the postprojective partition of a category ${\cal F}(\Delta)$ of $\Delta$-good modules which is finite. We shall see that this bound is given by the depths related to ${\cal F}(\Delta)$ of the morphisms $\pi(i): P(i) \rightarrow \Delta(i)$ for all $i=1,\cdots,n$. Similarly we shall find a bound to the number of levels of the preinjective partition of ${\cal F}(\Delta)$ from the depths related to ${\cal F}(\Delta)$ of the morphisms $\beta(i): \Delta(i) \rightarrow T(i)$, $i=1, \cdots, n$.

\begin{defn}
Let ${\cal F}(\Delta)$ be finite. Then there exists $n \geq 0$ such that ${\bf P_n}(\Delta) \neq \emptyset$ and ${\bf P_{n+1}}(\Delta)=\emptyset$. We define $p(\Delta)=n$. Dually, we define $q(\Delta)$ as the greatest level of a non-empty element of $\{{\bf I_n}(\Delta)\}_{n\geq 0}$.

\end{defn}

\begin{prop}
\label{boundpost}
If ${\cal F}(\Delta)$ is finite then $p(\Delta) \leq {\rm max}\{{\rm dp}_{\Delta}(\pi(i)):i=1, \cdots , n \}$.
\begin{proof}
Assume $p(\Delta)=m$ and take $M \in {\bf P_m}(\Delta)$. We have then $M= \Delta(j)$ for some $j \in \{1, \cdots, n \}$: otherwise, there would exist a proper epimorphism $M \rightarrow \Delta(j)$ which would lead to $\Delta(j) \in {\bf P_{m+1}}(\Delta)$ and hence ${\bf P_{m+1}}(\Delta) \neq \emptyset$. Consider the epimorphism $\pi(j): P(j) \rightarrow \Delta(j)$, with $\Delta(j) \in {\bf P_m}(\Delta)$. Since $\Delta(j) \in {\bf P_m}(\Delta)$, by Proposition \ref{propdan}, we have $\pi(j) \in {\rm rad}_{\Delta}^{m}(P(j), \Delta(j))$ which gives us $p(\Delta)=m \leq {\rm dp}_{\Delta}(\pi(j))$. The proposition follows.

\end{proof}
\end{prop}

\begin{prop}
\label{boundpre}
If ${\cal F}(\Delta)$ is finite then $q(\Delta) \leq {\rm max}\{{\rm dp}_{\Delta}(\beta(i)):i=1, \cdots , n \}$.
\begin{proof}

Assume $q(\Delta)=m$ and take $M \in {\bf I_m}(\Delta)$. Dualizing the arguments of Proposition \ref{boundpost} we have $M=\Delta(j)$ for some $j \in \{1, \cdots, n \}$. Consider the morphism $\beta(j): \Delta(j) \rightarrow T(j)$ with $\Delta(j) \in {\bf I_m}(\Delta)$. Since $\Delta(j) \in {\bf I_{m}}(\Delta)$ we obtain a monomorphism $g: \Delta(j) \rightarrow T$ such that $g \in {\rm rad}_{\Delta}^{m}(\Delta(j), T)$ from the chain ($\ast$) in the proof of Proposition \ref{propreview}. Following the line of the proof of $c) \Rightarrow a)$ in Theorem \ref{thmcharac} we have that there exists $h \in {\rm Hom}(T,T(j))$ such that $\beta(j)=hg$ and, hence, $\beta(j) \in {\rm rad}_{\Delta}^{m}(\Delta(j), T(j))$. Therefore, $q(\Delta)=m \leq {\rm dp}_{\Delta}(\beta(j))$.

\end{proof}

\end{prop}

\section*{Acknowledgements} This paper is part of the PhD's thesis of the author \cite{tese}. The author acknowledges financial suppport from CNPq during the preparation of this work. The author also would like to thank the referee for the valuable comments which helped to improve the manuscript. 

%\makeatletter \renewcommand{\@biblabel}[1]{\hfill#1.}\makeatother
%\newcommand{\bysame}{\leavevmode\hbox to3em{\hrulefill}\,}

%\section*{References}

\end{document}